\documentclass[a4paper,12pt]{amsart}

\usepackage{amsthm, amsmath, latexsym, amsfonts}
\usepackage[all]{xy}
\usepackage{amssymb}
\usepackage{mathrsfs}

\newtheorem{Lemma}{Lemma}[section]
\newtheorem{Proposition}[Lemma]{Proposition}
\newtheorem{Theorem}[Lemma]{Theorem}

\newtheorem{Definition}[Lemma]{Definition}

\theoremstyle{remark}
\newcommand{\mgnbar}{\overline{M}_{g,n}}
\newcommand{\mgn}{M_{g,n}}
\newcommand{\mg}{M_g}

\title{On the tautological ring of $\mgn$}

\author{Gilberto Bini and Claudio Fontanari}

\email{gilberto.bini@unimi.it} \curraddr{{\sc Dipartimento di Matematica \\ Universit\`a degli Studi di Milano \\ Via C. Saldini 50 \\ 20133 Milano \\ Italy.}}

\email{fontanar@science.unitn.it}\curraddr{{\sc Dipartimento di Matematica \\  Universit\`a degli Studi di Trento\\ Via Sommarive 14 \\ 38123 Trento \\ Italy.}}

\thanks{This research was partially supported by PRIN 2012 "Geometria delle variet\`a algebriche",
by FIRB 2012 "Moduli spaces and Applications", and by GNSAGA of INdAM (Italy).}



\begin{document}

\begin{abstract}
We state and check the analogue of Faber's conjectures for the tautological ring of the moduli spaces 
$\mgn$ of $n$-pointed smooth curves of genus $g$.  
\end{abstract}

\maketitle

\section{Introduction}

Let $\mgnbar$ be the moduli space of \emph{stable} curves of genus $g$ with $n$ marked points. 
We consider tautological classes on $\mgnbar$ by following
\cite{AC1}. Indeed, let $\overline{\pi}_{n+1}: \overline{M}_{g,n+1} \to \mgnbar$ be the morphism forgetting the last marked
point. Denote by $\sigma_1, \ldots, \sigma_n$ the canonical sections of $\overline{\pi}_{n+1}$, and by $D_1, \ldots, D_n$ the divisors
in $\overline{M}_{g,n+1}$ they define. We let $\omega_{\overline{\pi}_{n+1}}$ be the relative dualizing sheaf and we set
\begin{eqnarray*}
\psi_i &=& c_1(\sigma_i^*(\omega_{\overline{\pi}_{n+1}}))\\
K &=& c_1(\omega_{\overline{\pi}_{n+1}}(\sum D_i))\\
\kappa_i &=& \overline{\pi}_{n+1,*}(K^{i+1})\\
\tilde{\kappa}_i &=& \overline{\pi}_{n+1, *}(c_1(\omega_{\overline{\pi}_{n+1}})^{i+1}).
\end{eqnarray*}
\noindent
The classes $\kappa_i$ and $\tilde{\kappa}_i$ are related as
\begin{equation}\label{relation}
\kappa_i = \tilde{\kappa}_i+ \sum_{j=1}^n \psi_j^i
\end{equation}
by \cite{AC1}, equation (1.5). 

Let now $\mgn \subset \mgnbar$ be the moduli space of \emph{smooth} curves. 

\begin{Definition}\label{tautological}
The tautological ring $R^*(\mgn)$ is the subring of the rational Chow ring $A^*(\mgn)$ generated by the images 
of the classes $\psi_1, \ldots, \psi_n$ and $\tilde{\kappa}_i$ under the restriction map $A^*(\mgnbar) \to A^*(\mgn)$. 
\end{Definition}

The following statement is just the natural extension of classical Faber's conjectures (see \cite{Fab} 
for the original case $n=0$ and \cite{Yin}, sections 3.1 and 3.2, for the case $n = 1$).

\begin{Theorem}\label{main}

(i) The tautological ring $R^*(\mgn)$ satisfies $R^i(\mgn)=0$ for $i > g-1$ and $R^{g-1}(\mgn) = \mathbb{Q}^n$.

(ii) The ring $R^*(\mgn)$ is generated by $\tilde{\kappa}_{1}, \ldots, \tilde{\kappa}_{\lfloor g/3 \rfloor}$ and
by $\psi_1, \ldots, \psi_n$. There are no relations among these classes in codimension $\le \lfloor g/3 \rfloor$.

\end{Theorem}

We point out that for $n > 1$ the ring $R^*(\mgn)$ cannot be Gorenstein with socle in codimension $g-1$,
simply because in this case we have $R^{g-1}(\mgn) \ne \mathbb{Q}$ according to (i).

We remark that part (i) is already well-known to the experts. Indeed, the vanishing $R^i(\mgn)=0$ for $i > g-1$ 
is usually referred to as \emph{Getzler's conjecture} or \emph{Ionel's theorem} (see for instance 
\cite{GV}, \S 5.2 on pp. 29--30) and it turns out to be a direct consequence of Theorem $\star$ by Graber e Vakil 
(see \cite{GV}, Theorem 1.1). On the other hand, the second statement of (i) was conjectured by Dimitri Zvonkine 
in the unpublished note \cite{Zvo} and then proven by Alexandr Buryak in the online manuscript \cite{Bur}.   

We are grateful to Eduard Looijenga, Rahul Pandharipande and Sergey Shadrin for crucial remarks on a previous 
version of this note.

We work over the complex field $\mathbb{C}$.

\section{Proof of part (ii)}

From Definition \ref{tautological} and the vanishing $R^i(\mgn)=0$ for $i > g-1$ in Theorem \ref{main} (i), 
it follows that the tautological ring is additively generated by monomials
$$
\prod_{i=1}^n \psi_i^{e_i} \prod_{j \ge 0} \tilde{\kappa}_j^{f_j}
$$
satisfying $\sum_{i=1}^n e_i + \sum_{j \ge 0} j f_j < g$ (cf. \cite{FP}, Proposition 2).   

We also point out the following easy fact.

\begin{Lemma}\label{pullback}
Let $\pi_{n+1}: M_{g,n+1} \to \mgn$ be the morphism forgetting the last marked
point. Then $\pi_{n+1}^*(\tilde{\kappa_i})= \tilde{\kappa_i}$.
\end{Lemma}

\proof 
By \cite{AC1}, equation (1.10), we have $\overline{\pi}_{n+1}^*(\kappa_i) = \kappa_i - \psi_{n+1}^i$ 
and by \cite{AC2}, Lemma 3.1 (ii), we have $\overline{\pi}_{n+1}^*(\psi_i) = \psi_i - \Delta$,
where $\Delta$ denotes a boundary divisor class. Hence from (\ref{relation}) it follows that 
\begin{eqnarray*}
\overline{\pi}_{n+1}^*(\tilde{\kappa_i}) &=& \overline{\pi}_{n+1}^*(\kappa_i) - 
\sum_{j=1}^n \overline{\pi}_{n+1}^*(\psi_j^i) =\\
&=& \kappa_i - \psi_{n+1}^i - \sum_{j=1}^n \overline{\pi}_{n+1}^*(\psi_j)^i = \\
&=& \kappa_i - \psi_{n+1}^i - \sum_{j=1}^n (\psi_j - \Delta)^i.
\end{eqnarray*}
By restriction to the open part $M_{g, n+1}$ we get 
\begin{eqnarray*}
\pi_{n+1}^*(\tilde{\kappa_i}) &=& \kappa_i - \psi_{n+1}^i - \sum_{j=1}^n \psi_j^i = \\
&=& \kappa_i - \sum_{j=1}^{n+1} \psi_j^i = \tilde{\kappa_i}.
\end{eqnarray*}

\endproof

\begin{Proposition}
The ring $R^*(\mgn)$ is generated by $\tilde{\kappa}_{1}, \ldots, \tilde{\kappa}_{\lfloor g/3 \rfloor}$ 
and by $\psi_1, \ldots, \psi_n$.
\end{Proposition}

\proof 
As in the proof of \cite{Ion}, Theorem 1.5, it is enough to show that 
$$
\tilde{\kappa}_a = p(\tilde{\kappa}_1, \ldots, \tilde{\kappa}_{\lfloor g/3 \rfloor})
$$
for every $a \ge \lfloor g/3 \rfloor$, where $p$ denotes a polynomial. In order to do so, we consider the forgetful
map $\pi: \mgn \to \mg$ and via Lemma \ref{pullback} we reduce ourselves to the case of $\mg$, where 
$\tilde{\kappa}_a = \kappa_a$ by (\ref{relation}) and the proof of \cite{Ion}, Theorem 1.5, applies. 
Namely, 
\begin{eqnarray*}
\tilde{\kappa}_a &=& \pi^*(\tilde{\kappa}_a) = \pi^*(p(\kappa_1, \ldots, \kappa_{\lfloor g/3 \rfloor})) =\\
& & p(\pi^*(\tilde{\kappa}_1), \ldots, \pi^*(\tilde{\kappa}_{\lfloor g/3 \rfloor})) = 
p(\tilde{\kappa}_1, \ldots, \tilde{\kappa}_{\lfloor g/3 \rfloor}).
\end{eqnarray*}
\endproof

\begin{Proposition}
There are no relations among $\tilde{\kappa}_{1}, \ldots, \tilde{\kappa}_{\lfloor g/3 \rfloor}, \psi_1, \ldots, \psi_n$
in degree up to $\lfloor g/3 \rfloor$.
\end{Proposition}

\proof 
From the work of Boldsen \cite{Bol} it follows that 
$$
H^*(\mg, \mathbb{Q}) \cong \mathbb{Q}[\kappa_1, \ldots, \kappa_{\lfloor g/3 \rfloor}]
$$
in degree up to $\lfloor g/3 \rfloor$. On the other hand, by a theorem of Looijenga (\cite{Loo}, see also \cite{ACG}, 
\S 19, Theorem (5.5) on p. 684) we have
$$
H^*(\mgn, \mathbb{Q}) \cong H^*(\mg, \mathbb{Q})[\psi_1, \ldots, \psi_n]
$$
in the stable range. Hence we have 
\begin{eqnarray*}
A^*(\mgn) \otimes \mathbb{Q} &\cong& H^*(\mgn, \mathbb{Q}) \\
&\cong& \mathbb{Q}[\kappa_1, \ldots, \kappa_{\lfloor g/3 \rfloor}, 
\psi_1, \ldots, \psi_n] \\
&\cong& \mathbb{Q}[\tilde{\kappa}_1, \ldots, \tilde{\kappa}_{\lfloor g/3 \rfloor}, \psi_1, \ldots, \psi_n]
\end{eqnarray*}
in degree up to $\lfloor g/3 \rfloor$.
\endproof


\begin{thebibliography}{99}

\bibitem[AC1]{AC1} E. Arbarello and M. Cornalba: Combinatorial and algebro-geometric cohomology classes on the moduli
spaces of curves. J. Algebr. Geom. 5 (1996), 705--749.

\bibitem[AC2]{AC2} E. Arbarello and M. Cornalba: Calculating cohomology groups of moduli spaces of curves via 
algebraic geometry. Publ. Math. Inst. Hautes \'Etud. Sci. 88 (1998), 97--127.

\bibitem[ACG]{ACG} E. Arbarello, M. Cornalba and P. A. Griffiths: Geometry of algebraic curves. Volume 2. 
Grundlehren der Mathematischen Wissenschaften 268. Springer, Berlin 2011. 

\bibitem[Bol]{Bol} S. K. Boldsen: Improved homological stability for the mapping class group with integral 
or twisted coefficients. Math. Z. 270 (2012), 297--329.

\bibitem[Bur]{Bur} A. Buryak: Top tautological group of $\mgn$. Available online at 
\url{http://www.moebiuscontest.ru/files/2012/buryak.pdf}.

\bibitem[Fab]{Fab} C. Faber: A Conjectural Description of the Tautological Ring of the Moduli Space of Curves.
In: Moduli of curves and abelian varieties. The Dutch intercity seminar on moduli. Braunschweig: Vieweg. 
Aspects Math. E33 (1999), 109--129.

\bibitem[FP]{FP} C. Faber and R. Pandharipande: Relative maps and tautological classes. J. Eur. Math. Soc. 
(JEMS) 7 (2005), 13--49.

\bibitem[GV]{GV} T. Graber and R. Vakil: Relative virtual localization and vanishing of tautological classes
on moduli spaces of curves. Duke Math. J. 130 (2005), 1--37.

\bibitem[Ion]{Ion} E.-N. Ionel: Relations in the tautological ring of $\mathcal{M}_{g}$.
Duke Math. J. 129 (2005), 157--186.

\bibitem[Loo]{Loo} E. Looijenga: Stable cohomology of the mapping class group with symplectic coefficients 
and of the universal Abel-Jacobi map. J. Algebr. Geom. 5 (1996), 135--150. 

\bibitem[Yin]{Yin} Q. Yin: On the tautological rings of $\mathcal{M}_{g,1}$ and its universal Jacobian.
arXiv: 1206.3783 (2012).

\bibitem[Zvo]{Zvo} D. Zvonkine: A conjectural structure of the tautological ring of $\mgn$. Unpublished note
(2009).

\end{thebibliography}
\end{document}